\documentclass[12pt]{article}
\usepackage{amsthm}
\usepackage{amsmath}
\usepackage{amsfonts}
\usepackage{amssymb}
\usepackage{t-angles}
\allowdisplaybreaks

\newcommand{\sstyle}{\scriptstyle}

\newcommand{\id}{\mbox{id}}

\newcommand{\de}{\delta}

\newcommand{\ve}{\varepsilon}
\newcommand{\D}{\Delta}
\newcommand{\ot}{\otimes}

\newcommand{\di}{{}_{1}}
\newcommand{\dii}{{}_{2}}
\newcommand{\diii}{{}_{3}}
\newcommand{\diiii}{{}_{4}}

\newcommand{\doo}{{}_{(0)}}  %%\doo==\moo
\newcommand{\doi}{{}_{(-1)}}
\newcommand{\doii}{{}_{(-2)}}

\newcommand{\lmod}[1]{{}_{#1}\!{\cal M}}           % 左 模
\newcommand{\lcom}[1]{{}^{#1}\!{\cal M}}           % 左 余 模
\def\rbiprod{{\cdot\kern-.33em\triangleright\!\!\!<}}
\def\lbiprod{{>\!\!\!\triangleleft\kern-.33em\cdot\, }}

\def\lrbiprod{{\ \cdot\kern-.60em\triangleright\kern-.33em\triangleleft\kern-.33em\cdot\, }}

\def\lprod{{>\!\!\!\triangleleft\kern-.33em\ \, }}

\newcommand{\allowpagebreak}
\allowdisplaybreaks

\newtheorem{Theorem}{Theorem}[section]
\newtheorem{Proposition}[Theorem]{Proposition}

\newtheorem{Lemma}[Theorem]{Lemma}

\title{A Proof of  Radford's Biproduct Theorem by Using Braided Diagrams}
\author{Tao Zhang}
\date{}

\begin{document}

\footnotetext{Keywords: Hopf algebras, Biproduct Theorem, Braided Diagrams.}
\footnotetext{2010 Mathematics Subject Classification: 16S40, 16T05}

 \maketitle

 \setcounter{section}{0}

 \vskip0.1cm

{\bf Abstract.}\ \ We give a proof of Radford's Biproduct Theorem in
S. Montgomery's book [Hopf Algebras and Their Actions on Rings, CBMS
82, AMS,1993.] by using Majid's braided diagram method and Yu.
Bespalov and V. Lyubashenko's "t-angles.sty" package.

%{\bf Keywords:}\ \  Hopf algebras, Biproduct Theorem, Braided Diagrams.

\section{Introduction}
There is an important construction in the theory of quantum groups
and Hopf algebras, which is called biproduct by Radford in
\cite{Ra85} and bosonization by Majid in \cite{Ma94}. Through this
construction, one can get an ordinary Hopf algebra $B\star H$ from a
braided Hopf algebras $B$ and a Hopf algebras $H$, see Theorem \ref{maintheorem} below.
This construction can be used to give classification of pointed Hopf
algebras through braided Hopf algebras, see \cite{AS02}.

% The condition was first found by Radford in {\bf Theorem 1}
\begin{Theorem}\label{maintheorem}
Let $H$ be a bialgebra,  and $B$ is an algebra in $\lmod H$ and a
coalgebra in $\lcom H$. Then  $B\star H$ becomes a bialgebra $\iff$
$B$ is a coalgebra in $\lmod H$ , an algebra in $\lcom H$, $\ve_B$
is an algebra map, $\de (1_B)=1_B\ot 1_B$, and the following
identities hold:

(1) $\de_B(ab)=\sum a\di(a\dii\doi\cdot b\di)\ot a\dii\doo b\dii$,

(2) $\sum h\di b\doi\ot h\dii\cdot b\doo=\sum(h\di \cdot b)\doi
h\dii\ot(h\di\cdot b)\doo$.

\noindent If also $B$ has an antipode $S_B$ and $H$ is a Hopf
algebra with antipode  $S_H$, then $B\star H$ is a Hopf algebra with
antipode $S(b\star h)=\sum(1_B\star S_H(b_{-1}h))(S_B b_0 \star
1_H)=\sum (S_H(a\doi h))\di \cdot S_B(a\doo)\ot (S_H(a\doi h))\dii$.
\end{Theorem}

This theorem appeared in Remark 10.6.6 in Montgomery's book
\cite{Mo93}. After this remark, she said that "Although this
formulation may be more natural, and has the advantage of using $H$
and not $H^{cop}$, the computations are more tiresome because of the
notation".

We use braided diagrams to overcome this difficulty. This method was
used by Majid in \cite[Theorem 2.4]{Ma94}, but in that paper he
assume $H$ to be cocommutative. We find that this condition is not
necessary and the same method can also be applied.

Throughout the paper we freely use the notations and conventions of
\cite{Ra85, Mo93} with slightly differences. In particular all
vector spaces will be over a field $k$. For an algebra $A$ with
multiplication $m: A\ot A\to A$, we write $m(a\otimes b)=ab$ for simplicity.
For a coalgebra $C$ with
comultiplication $\D:C\to  C\ot C$, we write 
$\D(c)=\sum c\di\ot c\dii$ using Sweedler's notation.
\section{Preliminaries}

Let $(H,m)$ be an associative algebra and $(H,\D)$ be a coassociative coalgebra over field $k$. Then $H$ is called a bialgebra  if the following compatible 
conditions hold:
\begin{equation}\label{campatible}
\D(hg)=\D(h)\D(g), \quad \D(1_H)=1_H\ot 1_H,\quad \varepsilon(1_H)=1_H.
\end{equation}
It is called a Hopf algebra if there is an antipode $S: H\to H$ such that 
\begin{equation}\label{antipode}
(\id\ot S)\Delta(h)=\varepsilon(h) 1_H=(S\ot \id)\Delta(h).
\end{equation}

A left $H$-module is a $k$-vector space $M$ with a $k$-linear map $\alpha: H\ot M\to M:
\alpha(h\ot m)=h\cdot m$ such that $(gh)\cdot m=h\cdot (g\cdot m)$.
A left $H$-comodule is a $k$-vector space $M$ with a $k$-linear map $\rho:
M\to H\ot M: \rho(m)=\sum m\doi\ot m\doo$ such that
\begin{equation}\label{com}
\sum m\doi\ot m\doo\doi\ot m\doo\doo=\sum m\doi\di\ot m\doi\dii \ot
m\doo.\end{equation}

Here we review some basic facts about algebras and coalgebras in the
category $\lmod H$ of left $H$-modules and in the category $\lcom H$
of left $H$-comodules. First of all recall that if $M$ and $N$ are
left $H$-modules, then the left $H\ot H$-modules $M\ot N$ is also a
left $H$-modules by pull-back along $\D$, i.e., $h\cdot(m\ot n)=\sum
h\di\cdot m\ot h\dii\cdot n$.

An algebra $B$ in $\lmod H$ is an left $H$-module algebra, that is,
$B$ is a left $H$-module and also a $k$-algebra $(B, m, \eta)$ such
that $m$ and $\eta$ are module maps, i.e.
\begin{equation}\label{ma}
h\cdot (ab)=\sum (h\di \cdot a)(h\dii \cdot b),\quad h\cdot 1_B=\ve_H(h)1_B,\quad \forall h\in H, a,b\in B.
\end{equation}

 An coalgebra $B$ in $\lmod H$ is an $H$-module coalgebra,
that is, $B$ is a left $H$-module and also a $k$-coalgebra $(B, \D,
\ve )$ such that $\D$ and $\ve $ are module maps, i.e.
\begin{equation}\label{mc}
\D(h\cdot b)=\sum (h\di \cdot b\di)\ot(h\dii\cdot b\dii),\quad \ve_B(h\cdot b)=\ve_H(h)\ve_B(b), \quad \forall h\in H, b\in B.
\end{equation}

An algebra $B$ in $\lcom H$ is an $H$-comodule algebra, that is, $B$
is a left $H$-comodule and also a $k$-algebra $(B, m, \eta)$ such
that $m$ and $\eta$ are module maps, i.e.
\begin{equation}\label{ca}
\rho(ab)=\sum a\doi b\doi\ot a\doo b\doo,\quad \rho(1_H)=1_H\ot 1_B,\quad\forall a,b\in B.
\end{equation}

An coalgebra $B$ in $\lcom H$ is an $H$-comodule coalgebra, that is,
$B$ is a left $H$-comodule and also a  $k$-coalgebra $(B, \D, \ve )$
such that $\D$ and $\ve $ are comodule maps, i.e.
\begin{equation}\label{cc}
\sum b\doi\ot b\doo\di\ot b\doo\dii=\sum b\di\doi b\dii\doi\ot
b\di\doo\ot b\dii\doo, \quad \sum\ve_B(b\doo)b\doi=\ve_B(b)1_H,\quad \forall b\in B.
\end{equation}

\section{The Proof}

In this section, we give a detailed proof of Theorem
\ref{maintheorem}. First, we review some facts about the product and
coproduct in $B\star H$.

Let $B$ be an left $H$-module algebra. Then the smash product
algebra $B\# H$ is defined as follows: $B\# H:=B\ot H$ as a vector space,
with multiplication given by
$$(a\# g)(b\# h)=\sum a(h\di \cdot b)\# h\dii g,$$
and unit $1_{B\# H}=1_B\ot 1_H,$ for all $a, b\in B, h,g\in H$.

Let $H$ be a bialgebra and $B$ a coalgebra in $\lcom H$. The smash
coproduct $B\sharp H$ is defined to be $B\ot H$ as a vector space,
with comultiplication given by
$$\D(b\sharp h)=\sum b\di\sharp b\dii\doi h\di\ot b\di\doo\sharp h\dii,$$
and counit $\ve_{B\sharp H} (b\sharp h)=\ve_B (b)\ve_H (h),$ for all $b\in B, h\in
H$.

\begin{Lemma}\label{lemmma} $B\sharp H$ is a coalgebra with the above
comultiplication.
\end{Lemma}
\proof We check that $\D$ is coassociative. Now for all $b\in B, h\in
H$,
\begin{eqnarray*}
&&(\Delta\ot \id) \Delta(b\sharp h))\\
&=&\sum b\di\di \sharp b\di\dii\doi (b\dii\doi h\di)\di\ot b\di\dii\doo \sharp (b\dii\doi h\di)\dii\ot b\dii\doo\sharp h\dii \\
&=&\sum b\di \sharp b\dii\doi b\diii\doii h\di\ot b\dii\doo \sharp b\diii\doi h\dii\ot b\diii\doo\sharp h\diii \\
&& \mbox{by equation (\ref{com}) applied to $b_3$}\\
&=&\sum b\di \sharp b\dii\doi b\diii\doi h\di\ot b\dii\doo \sharp b\diii\doo\doi h\dii\ot b\diii\doo\doo\sharp h\diii \\
&& \mbox{ by equation (\ref{cc}) applied to $b_2$}\\
&=&\sum b\di \sharp b\dii\doi  h\di\ot b\dii\doo \sharp b\dii\doo\dii\doi h\dii\ot b\dii\doo\dii\doo\sharp h\diii \\
&=&(\id\ot \Delta) \Delta(b\sharp h)).
\end{eqnarray*}
Thus we get that the comultiplication is coassociative. \qed

{\bf Proof of Theorem \ref{maintheorem} } We will show that $B\star
H$ is a bialgebra. First, $B\star H$ is an algebra by smash product,
and it is a coalgebra by Lemma \ref{lemmma}. We now check the compatible condition:
%$\D$ is multiplicative.
\begin{eqnarray*}
&&\Delta((a\star g)(b\star h))=\Delta(a(g\di b) \star g\dii h)\\
&=&\sum (a(g\di b))\di \star (a(g\di b))\dii\doi (g\di h)\di\ot(a(g\di b))\dii\doo \star (g\dii h)\dii \\
&& \mbox{by condition (1) in Theorem \ref{maintheorem}  applied to $a(g\di b)$}\\
&=&\sum a\di(a\dii\doi(g\di b)\di) \star (a\dii\doi(g\di b)\dii)\doi(g\dii h)\di\ot (a\dii\doi(g\di b)\dii)\doo \star g\dii\dii h\dii \\
&& \mbox{by equation (\ref{mc}) applied to $g\di \cdot b$ and equation (\ref{ca}) applied to $a\dii\doo (g\di b)\dii$}\\
&=&\sum a\di(a\dii\doi(g\di\di b\di)) \star (a\dii\doi(g\di\dii b\dii))\doi(g\dii\di h\di)\ot (a\dii\doi(g\di\dii b\dii))\doo \star g\dii\dii h\dii\\
&=&\sum a\di(a\dii\doi(g\di b\di)) \star a\dii\doo\doi (g\dii
b\dii)g\diii h\di\ot
a\dii\doo\doo(g\dii h\dii)\doo \star g\diiii h\dii \\
&=&\sum a\di(a\dii\doi(g\di b\di)) \star a\dii\doo\doi g\dii
b\dii\doi
h\di\ot a\dii\doo\doo(g\diii h\dii\doo) \star g\diiii h\dii \\
&& \mbox{by condition (2) in Theorem \ref{maintheorem}  applied to $g\dii\ot b\dii$}\\
&=&\sum a\di((a\dii\doi\di g\di b\di) \star a\dii\doi\dii g\dii
b\dii\doi
h\di\ot a\dii\doo(g\diii h\dii\doo) \star g\diiii h\dii\\
&=&\sum a\di((a\dii\doi g)\di b\di) \star (a\dii\doi g\di)\dii
b\dii\doi h\di\ot a\dii\doo(g\diii h\dii\doo) \star g\diiii h\dii\\
&=&\Delta(a\star g)\Delta(b\star h).
\end{eqnarray*}
%%%%%%%%%%%%%%%%%%%%%%%%%%%%%%%%%%%%%%%%%%%%%%%%%%%%%%%%%%%%%%%
The fact that the given $S$ make $B\star H$ into a Hopf algebra is checked below:
\begin{eqnarray*}
&&(S\ot\id) \D(b\star h)=S(b\star h)\di (b\star h)\dii\\
&=&(S_H(b\di\doi b\dii \doi h\di)\di\cdot
S_B(b\di\doo))(S_H(b\di\doi b\dii \doi h\di)\dii\cdot b\dii\doo)\ot S_H(b\di\doi b\dii \doi h\di)\diii h\dii\\
%&& \mbox{[by definition of $B\star H$]}\\
&=&(S_H(b\doi h\di)\di\cdot S_B(b\doo\di)(S_H(b\doi h\di)\dii\cdot
b\doo\dii)\ot S_H(b\di\doi h\di)\diii h\dii\\
&& \mbox{by $B$ is an $H$-comodule coalgebra}\\
&=&(S_H(b\doi h\di)\di\cdot (S_B(b\doo\di)b\doo\dii)\ot S_H(b\doi h\di)\dii h\dii\\
&& \mbox{by $B$ is an $H$-module algebra applied to $S_H(b\doi h\di)\di \ot S_B(b\doo\di)\ot b\doo\doi$}\\
&=& \ve_B (b\doo)(S_H(b\doi h\di)\di\cdot 1_H)\ot S_H(b\doi h\di)\dii h\dii\\
&& \mbox{by antipode of $B$ applied to $b\doo$}\\
&=& \ve_B (b)\ve_B((S_H(1_H h\di))\di)1_B\ot S_H(1_H h\di)\dii h\dii\\
% && \mbox{[by  ********  applied to $b\doo$]}\\
&=& \ve_B (b)\ve_B(S_H(1_H h\di)\di 1_H\ot S_H(1_H h\di)\dii h\dii\\
&=& \ve_B (b) 1_B\ot S(h\di)h\dii\\
% && \mbox{[by antipode of $H$ applied to $h$]}\\
&=& \ve_B(b)\ve (h)  1_B\ot  1_H.\\
%\end{eqnarray*}
%\begin{eqnarray*}
&&\\[-1em]
&&(\id \ot S) \D(b\star h)=(b\star h)\di S(b\star h)\dii\\
&=& b\di ((b\dii\doi h\di)\di \cdot(S_H(a\dii\doo\doi h\dii)\di\cdot
S_B(a\dii\doo\doo))\ot (b\dii\doi h\di)\dii S_H(a\dii\doo\doi
h\dii)\dii\\
&=& b\di ((b\dii\doi \di h\di)\di S_H(b\dii\doi\di h\dii)\di\cdot S_B(b\dii\doo))\ot (b\dii\doi \di h\di)\doi S_H(b\dii\doi\di h\dii)\dii\\
&=& b\di (((b\dii\doi \di h)\di S_H(b\dii\doi\di h)\dii)\di\cdot S_B(b\dii\doo))\ot (b\dii\doi \di h)\di S_H(b\dii\doi\di h)\dii)\dii\\
&=&\ve_B (b\dii\doi h)b\di(1_H\cdot S_B(b\dii\doo)) \ot  1_H\\
&=&\ve_B (b\dii\doi) b\di (1_H\cdot S(b\dii\doo))\ot \ve_H (h) 1_H \\
&=&b\di S(b\dii)\ot \ve_H (h) 1_H \\
% && \mbox{[by antipode of $H$ applied to $h$]}\\
&=& \ve_B(b)\ve (h)  1_B\ot  1_H.
\end{eqnarray*}
The proof is completed. \qed

%we give details of all proofs for the reader's convenience.
When we apply the biproduct construction to the category of
quasitiangular Hopf algebras, we get the following result. It is due
to Majid \cite{Ma94} which is called "bosonization".
\begin{Proposition}\label{prop1}
If $(H, R)$ be a quasitriangular bialgebra and let $B$ be a
bialgebra in $\lmod H$. Then $B$ is also a left $H$-comodule algebra
by defining $\rho: B\to H\ot B $ via $\rho(b)=R^{-1}(1\ot b)$ for
all $b\in B$. Then the biproduct $B\star H$ is a bialgebra. If also
$H$ is a Hopf algebra and $B$ is a Hopf algebra in $\lmod H$, then
$B\star H$ is a Hopf algebra.
\end{Proposition}

%\begin{Proposition}\label{prop2}
%If $x_B\in B$ and $x_H\in H$ are right integrals, then $x_B\star x_H$ is  a right integral of $B\star H$.
%\end{Proposition}

\section{Braided diagrams: the method we use}
In this section, we give an introduction to the method of braided
diagrams in Hopf algebras from which we abtained the above results. This
method was introduced by D. N. Yetter in \cite{Ye90} where he found
that the category of Yetter-Drinfeld modules form a braided monoidal
category. There he call it crossed bimodules which  in fact is
condition (2) of Theorem \ref{maintheorem}. Braided diagram was used by many
authors such as S.Majid \cite{Ma94}, Y.Bespalov and B. Drabant
\cite{BD99}, S.C. Zhang and H.X. Chen \cite{ZC01}, etc. Here we
use Yu. Bespalov and V. Lyubashenko's "t-angles.sty" package as in
\cite{BD99} to draw the diagrams.

First, all maps are written downwards from top to bottom. The maps
$m$, $\D$, $\alpha$, $\rho$, $\ve$, $\eta$ are graphically written as
\[\hstretch 50 \vstretch 50
m_H=
\begin{tangle}
\object{\sstyle H}\step[2]\object{\sstyle H}\\
\cu\\
\step\object{\sstyle H}\\
\end{tangle}
\enspace, \step \delta_H=
\begin{tangle}
\step\object{\sstyle H}\\
\cd\\
\object{\sstyle H}\step[2]\object{\sstyle H}\\
\end{tangle}
\enspace, \step m_B=
\begin{tangle}
\object{\sstyle B}\step[2]\object{\sstyle B}\\
\cu\\
\step\object{\sstyle B}\\
\end{tangle}
\enspace, \step \delta_B=
\begin{tangle}
\step\object{\sstyle B}\\
\cd\\
\object{\sstyle B}\step[2]\object{\sstyle B}\\
\end{tangle}
\enspace, \step\alpha=
\begin{tangle}
\object{\sstyle H}\step[2]\object{\sstyle B}\\
\Lu\\
\step[2]\object{\sstyle B}\\
\end{tangle}
\enspace, \step \rho=
\begin{tangle}
\step[2]\object{\sstyle B}\\
\Ld\\
\object{\sstyle H}\step[2]\object{\sstyle B}\\
\end{tangle}
\enspace,
\step \ve_H=
\begin{tangle}
\object{\sstyle H}\\
\id\\
\QQ \epsilon \\
\end{tangle}
\enspace,
\step \eta_H=
\begin{tangle}
\Q \eta\\
\id\\
\object{\sstyle H}\\
\end{tangle}
\]
%%%%%%%%%%%%%%%%%%%%%%%%%%%%%%%%%%%%%%%%%

Secondly, we can picture the conditions of module algebra, module
coalgebra, comodule algebra, comodule coalgebra, algebra-coalgebra,
Yetter-Drinfel'd module as follows:
\[
\hstretch 50 \vstretch 50
\begin{tangle}
\object{\sstyle H}\step[1]\object{\sstyle B}\step[2]\object{\sstyle B}\\
\id\step[1]\cu\\
\Lu\\
\step[2]\object{\sstyle B}
\end{tangle}
\enspace=\enspace
\begin{tangle}
\step\object{\sstyle H}\step[3]\object{\sstyle B}\step[2]\object{\sstyle B}\\
\cd\step[2]\id\step[2]\id\\
\id\step[2]\x\step[2]\id\\
\Lu\step[2]\Lu\\
\step[2]\Cu\\
\step[4]\object{\sstyle B}
\end{tangle}\enspace,\step[2]
\begin{tangle}
\object{\sstyle H}\step[2]\object{\sstyle B}\\
\Lu\\
\step\cd\\
\step\object{\sstyle B}\step[2]\object{\sstyle B}
\end{tangle}
\enspace=\enspace
\begin{tangle}
\step\object{\sstyle H}\step[3]\object{\sstyle B}\\
\cd\step\cd\\
\id\step[2]\hx\step[2]\id\\
\Lu\step\Lu\\
\Step\object{\sstyle B}\step[3]\object{\sstyle B}
\end{tangle}\enspace,\step[2]
\begin{tangle}
\step\object{\sstyle B}\step[2]\object{\sstyle B}\\
\step\cu\\
\Ld\\
\object{\sstyle H}\step[2]\object{\sstyle B}
\end{tangle}
\enspace=\enspace
\begin{tangle}
\step[2]\object{\sstyle B}\step[3]\object{\sstyle B}\\
\Ld\step\Ld\\
\id\step[2]\hx\step[2]\id\\
\cu\step\cu\\
\step\object{\sstyle H}\step[3]\object{\sstyle B}
\end{tangle}
\enspace,\step[2]
\begin{tangle}
\step[2]\object{\sstyle B}\\
\Ld\\
\id\step[1]\cd\\
\object{\sstyle H}\step[1]\object{\sstyle B}\step[2]\object{\sstyle
B}
\end{tangle}
\enspace=\enspace
\begin{tangle}
\step[4]\object{\sstyle B}\\
\step[2]\Cd\\
\Ld\step[2]\Ld\\
\id\step[2]\x\step[2]\id\\
\cu\step[2]\id\step[2]\id\\
\step\object{\sstyle H}\step[3]\object{\sstyle
B}\step[2]\object{\sstyle B}
\end{tangle}
\enspace.
\]
%%%%%%%%%%%%%%%%%%%%%%%%%%%%%%%%%%%%%%
\[
\hstretch 50 \vstretch 50
\begin{tangle}
\object{\sstyle H}\step[2]\object{\sstyle H}\\
\cu\\
\cd\\
\object{\sstyle H}\step[2]\object{\sstyle H}
\end{tangle}
\enspace=\enspace
\begin{tangle}
\step\object{\sstyle H}\step[3]\object{\sstyle H}\\
\cd\step\cd\\
\id\step\step\hx\step[2]\id\\
\cu\step\cu\\
\step\object{\sstyle H}\step[3]\object{\sstyle H}
\end{tangle}
\enspace,\step[2]
\begin{tangle}
\object{\sstyle B}\step[2]\object{\sstyle B}\\
\cu\\
\cd\\
\object{\sstyle B}\step[2]\object{\sstyle B}
\end{tangle}
\enspace=\enspace
\begin{tangle}
\step\object{\sstyle B}\step[3]\object{\sstyle B}\\
\cd\step\cd\\
\id\step\ld\step\id\step[2]\id\\
\id\step\id\step\hx\step[2]\id\\
\id\step\lu\step\id\step[2]\id\\
\cu\step\cu\\
\step\object{\sstyle B}\step[3]\object{\sstyle B}
\end{tangle}
\enspace,\step[2]
\begin{tangle}
\step\object{\sstyle H}\step[4]\object{\sstyle B}\\
\cd\step\Ld\\
\id\step[2]\hx\step[2]\id\\
\cu\step\Lu\\
\step\object{\sstyle H}\step[4]\object{\sstyle B}
\end{tangle}
\enspace=\enspace
\begin{tangle}
\step\object{\sstyle H}\step[3]\object{\sstyle B}\\
\cd\step[2]\id\\
\x\step[2]\id\\
\id\step[2]\Lu\\
\id\step[2]\Ld\\
\x\step[2]\id\\
\cu\step[2]\id\\
\step\object{\sstyle H}\step[3]\object{\sstyle B}
\end{tangle}
\enspace=\enspace
\begin{tangle}
\step\object{\sstyle H}\step[3]\object{\sstyle B}\\
\cd\step[2]\id\\
\id\step[2]\x\\
\Lu\step[2]\id\\
\Ld\step[2]\id\\
\id\step[2]\x\\
\cu\step[2]\id\\
\step\object{\sstyle H}\step[3]\object{\sstyle B}
\end{tangle}
\enspace.
\]
\[
\hstretch 50 \vstretch 50 m=\
\begin{tangle}
\object{\sstyle B}\step[2]\step\object{\sstyle H}\step\step[2]\object{\sstyle B}\step[2]\object{\sstyle H}\\
\id\step[2]\cd\Step\id\Step\id\\
\id\step[2]\id\Step\x\Step\id\\
\id\step[2]\lu[2] \step[2]\id\step[2]\id\\
\Cu\Step\cu\\
\step[2]\object{\sstyle B}\step[4]\step\object{\sstyle H}\\
\end{tangle}
\enspace,\step[2]\D=\
\begin{tangle}
\step[2]\object{\sstyle B}\step[4]\step\object{\sstyle H}\\
\Cd\Step\cd\\
\id\step[2]\ld[2] \step[2]\id\step[2]\id\\
\id\step[2]\id\Step\x\Step\id\\
\id\step[2]\cu\Step\id\Step\id\\
\object{\sstyle B}\step[2]\step\object{\sstyle H}\step\step[2]\object{\sstyle B}\step[2]\object{\sstyle H}\\
\end{tangle}
\enspace,\step[2] S=\
\begin{tangle}
\step[2]\object{\sstyle B}\step[2]\object{\sstyle H}\\
\Ld\Step\id\\
\id\step[2]\x\\
\cu\Step\id\\
\obox 2{\sstyle {S_H}}\step\obox 2{\sstyle {S_B}}\\
\cd\Step\id\\
\id\step[2]\x\\
\Lu\Step\id\\
\step[2]\object{\sstyle B}\step[2]\object{\sstyle H}\\
\end{tangle}\ .
\]

Finally, we give the diagrammatic proof of Lemma \ref{lemmma} and
Theorem \ref{maintheorem}.
\begin{figure}[!ht]
\[
\hstretch 50 \vstretch 50
\begin{tangle}
\step[6]\object{\sstyle B}\step[4]\step\object{\sstyle H}\\
\step[4]\Cd\Step\cd\\
\step[3]\dd\step[2]\ld[2] \step[2]\id\step[2]\id\\
\step[2]\dd\step[3]\id\Step\x\Step\id\\
\step[2]\id\step[4]\cu\Step\id\Step\id\\
\Cd\Step\cd\step[2]\id\step[2]\id\\
\id\step[2]\ld[2] \step[2]\id\step[2]\id\step[2]\id\step[2]\id\\
\id\step[2]\id\Step\x\Step\id\step[2]\id\step[2]\id\\
\id\step[2]\cu\Step\id\Step\id\step[2]\id\step[2]\id\\
\object{\sstyle B}\step[2]\step\object{\sstyle H}\step\step[2]
\object{\sstyle B}\step[2]\object{\sstyle H}\step[2]\object{\sstyle B}\step[2]\object{\sstyle H}\\
\end{tangle}
\enspace=\step[.5]
\begin{tangle}
\step[6]\object{\sstyle B}\step[4]\step\object{\sstyle H}\\
\step[4]\Cd\Step\cd\\
\step[3]\dd\step[2]\ld[2] \step[2]\id\step[2]\id\\
\step[2]\dd\step[2]\dd\Step\x\Step\id\\
\step[2]\id\step[2]\cd\step\cd\step\id\Step\id\\
\step[2]\id\step[2]\id\Step\hx\Step\id\step\id\step\step\id\\
\Cd\cu\step\cu\step\id\step[2]\id\\
\id\step[2]\ld[2] \step[1]\d\step[2]\id\step[2]\id\step[2]\id\\
\id\step[2]\id\Step\x\Step\id\step[2]\id\step[2]\id\\
\id\step[2]\cu\Step\id\Step\id\step[2]\id\step[2]\id\\
\object{\sstyle B}\step[2]\step\object{\sstyle H}\step\step[2]
\object{\sstyle B}\step[2]\object{\sstyle H}\step[2]\object{\sstyle B}\step[2]\object{\sstyle H}\\
\end{tangle}
\enspace=\enspace
\begin{tangle}
\step[3]\object{\sstyle B}\step[4]\step\object{\sstyle H}\\
\step[1]\Cd\step[3]\id\\
\cd\step[1]\ld[2] \step[2]\cd\\
\step[-1]\dd\step[1]\ld\step\id\step\ld\step[1]\cd\step[1]\id\\
\step[-1]\id\step[1]\dd\step\hx\step[1]\id\step[1]\hx\step[2]\id\step[1]\id\\
\step[-1]\id\step[1]\id\step\dd\step\id\step\hx\step[2]\se2\step[-1]\dd\step[1]\id\\
\step[-1]\id\step[1]\id\step\id\step[2]\hx\step\cu\step\id\step\id\\
\step[-1]\id\step[1]\id\step[1]\cu\step\id\step[2]\id\step\step\id\step\id\\
\step[-1]\id\step[1]\cu\step[2]\id\step[2]\id\step[1]\step\id\step\id\\
\step[-1]\object{\sstyle B}\step[1]\step\object{\sstyle H}
\step\step[2]\object{\sstyle B}\step[2]\object{\sstyle H}\step[2]\object{\sstyle B}\step[1]\object{\sstyle H}\\
\end{tangle}
\enspace=\enspace
\begin{tangle}
\step[2]\object{\sstyle B}\step[6]\step\object{\sstyle H}\\
\Cd\step[5]\id\\
\id\step[3]\cd \step[4]\id\\
\id\step[2]\ld\step\ld\step[3]\cd\\
\id\step[1]\dd\step\hx\step[1]\id\step[2]\dd\step[2]\d\\
\id\step[1]\cu\step\id\step[2]\ne2\step[-2]\nw3\step[5]\id\\
\id\step[2]\id\step[2]\hx\step\ld[2]\step[2]\cd\\
\id\step[2]\cu\step\id\step\id\step[2]\x \step[2]\id\\
\id\step[2]\step[1]\id\step[2]\id\step[1]\cu\Step\id\Step\id\\
\object{\sstyle B}\step[2]\step\object{\sstyle H}\step\step[2]
\object{\sstyle B}\step[2]\object{\sstyle H}\step[2]\object{\sstyle B}\step[2]\object{\sstyle H}\\
\end{tangle}
\enspace=\enspace
\begin{tangle}
\step[2]\object{\sstyle B}\step[4]\step\object{\sstyle H}\\
\Cd\Step\cd\\
\id\step[2]\ld[2] \step[2]\id\step[2]\d\\
\id\step[2]\id\Step\x\step[3]\d\\
\id\step[2]\cu\Step\id\step[4]\d\\
\id\step[2]\dd\step[1]\Cd\Step\cd\\
\id\step[2]\id\step[2]\id\step[2]\ld[2] \step[2]\id\step[2]\id\\
\id\step[2]\id\step[2]\id\step[2]\id\Step\x\Step\id\\
\id\step[2]\id\step[2]\id\step[2]\cu\Step\id\Step\id\\
\object{\sstyle B}\step[2]\step\object{\sstyle H}\step\step[2]
\object{\sstyle B}\step[2]\object{\sstyle H}\step[2]\object{\sstyle B}\step[2]\object{\sstyle H}\\
\end{tangle}
\]
\caption{Proof of coalgebra structure of $B\star H$}
\label{fig:coproduct}
\end{figure}
%**********************************************************************
\begin{figure}[!ht]
\[
\hstretch 50 \vstretch 50
\begin{tangle}
\object{\sstyle B}\step[2]\step\object{\sstyle H}\step\step[2]\object{\sstyle B}\step[2]\object{\sstyle H}\\
\id\step[2]\cd\Step\id\Step\id\\
\id\step[2]\id\Step\x\Step\id\\
\id\step[2]\lu[2] \step[2]\id\step[2]\id\\
\Cu\Step\cu\\
\Cd\Step\cd\\
\id\step[2]\ld[2] \step[2]\id\step[2]\id\\
\id\step[2]\id\Step\x\Step\id\\
\id\step[2]\cu\Step\id\Step\id\\
\object{\sstyle B}\step[2]\step\object{\sstyle H}\step\step[2]\object{\sstyle B}\step[2]\object{\sstyle H}\\
\end{tangle}
\hstretch 100 \vstretch 100
\enspace=\enspace
\begin{tangle}
\step[.5]\object{\sstyle B}\step[1.5]\object{\sstyle H}\step[1.5]\object{\sstyle B}\step[2]\object{\sstyle H}\\
\hh\hstep\id\step\cd\step\id\step[2]\id\\
\hh\hstep\id\step\id\step\x\step[2]\id\\
\hh\cd\hstep\lu\step\d\step[1.5]\id\\
\hh\id\hstep\hld\step\cd\hstep\cd\hstep\cd\\
\hh\id\hstep\id\hstep\x\step\id\hstep\id\step\d\step[-.5]\dd\step\id\\
\hh\id\hstep\hlu\step\cu\hstep\cu\hstep\cu\\
\hh\cu\hstep\ld\step\dd\step[1.5]\id\\
\hh\hstep\id\step\id\step\x\step[2]\id\\
\hh\hstep\id\step\cu\step\id\step[2]\id\\
\step[.5]\object{\sstyle B}\step[1.5]\object{\sstyle H}\step[1.5]\object{\sstyle B}\step[2]\object{\sstyle H}
\end{tangle}
\enspace=\enspace
\begin{tangle}
\object{\sstyle B}\step[2]\object{\sstyle H}\step[1.5]\object{\sstyle B}\step[2.5]\object{\sstyle H}\\
\hh\hstep\step[-.5]\id\hstep\step\cd\step\id\step[2.5]\id\\
\hh\hstep\step[-.5]\id\hstep\hstep\hstep\id\step\se1\step[-.5]\dd\step[2.5]\id\\
\hh\step[-.5]\cd\step[.5]\cd\hstep\cd\se1\step[1.5]\id\\
\hh\step[-.5]\id\hstep\hld\step[.5]\id\hstep\sw1\step[-.5]\nw1\hstep\id\step\id\step\hstep\id\\
\hh\step[-.5]\id\step[.5]\id\hstep\id\hstep\hlu\step\hstep\hlu\hstep\cd\hstep\cd\\
\hh\step[-.5]\id\hstep\id\hstep\x\step[.5]\hstep\step\id\hstep\id\step\d\step[-.5]\dd\hstep\hstep\id\\
\hh\step[-.5]\id\step[.5]\id\hstep\id\hstep\hld\step\hstep\hld\hstep\cu\hstep\cu\\
\hh\step[-.5]\id\hstep\hlu\step[.5]\id\hstep\nw1\step[-.5]\sw1\hstep\id\step\id\step\hstep\id\\
\hh\step[-.5]\cu\step[.5]\cu\hstep\cu\ne1\step[1.5]\id\\
\hh\hstep\step[-.5]\id\hstep\hstep\hstep\id\step\ne1\step[-.5]\d\step[2.5]\id\\
\hh\hstep\step[-.5]\id\hstep\step\cu\step\id\step[2.5]\id\\
\object{\sstyle B}\step[2]\object{\sstyle H}\step[1.5]\object{\sstyle B}\step[2.5]\object{\sstyle H}\\
\end{tangle}
\enspace=\enspace \hstretch 50 \vstretch 50
\begin{tangle}
\object{\sstyle B}\step[4]\object{\sstyle H}\step[4]\object{\sstyle B}\step[4]\object{\sstyle H}\\
\id\step[2]\Cd\step[2]\id\step[4]\id\\
\id\step[1]\cd\step[3]\x\step[4]\id\\
\id\step[1]\id\step[2]\id\step[2]\cd\step[1]\d\step[3]\id\\
\id\step[1]\id\step[2]\x\step[2]\d\step[1]\d\step[2]\id\\
\step[-1]\cd\Lu\step[1]\cd\step[2]\id\step[2]\id\step[2]\id\\
\step[-1]\id\step\ld\step[2]\id\step[1]\id\step[2]\x\step[2]\id\step[2]\id\\
\step[-1]\id\step\id\step[1]\id\step[2]\id\step[1]\Lu\step[2]\id\step[2]\id\step[1]\cd\\
\step[-1]\id\step\id\step\x\step[1]\Ld\step[2]\id\step[2]\hx\step[2]\id\\
\step[-1]\id\step\lu\step[2]\id\step[1]\id\step[2]\x\step[2]\id\step[1]\cu\\
\step[-1]\cu\Ld\step[1]\cu\step[2]\id\step[2]\id\step[2]\id\\
\id\step[1]\id\step[2]\x\step[2]\dd\step[1]\dd\step[2]\id\\
\id\step[1]\id\step[2]\id\step[2]\cu\step[1]\dd\step[3]\id\\
\id\step[1]\cu\step[3]\x\step[4]\id\\
\id\step[2]\Cu\step[2]\id\step[4]\id\\
\object{\sstyle B}\step[4]\object{\sstyle H}\step[4]\object{\sstyle B}\step[4]\object{\sstyle H}\\
\end{tangle}
\]
\[
\enspace=\enspace
\begin{tangle}
\object{\sstyle B}\step[2]\object{\sstyle H}\step[1.5]\object{\sstyle B}\step[2]\object{\sstyle H}\\
\hh\hstep\step[-.5]\id\hstep\step\cd\step\id\step[2]\id\\
\hh\hstep\step[-.5]\id\hstep\hstep\hstep\id\step\se1\step[-.5]\dd\step[2]\id\\
\hh\step[-.5]\cd\step[.5]\cd\hstep\cd\se1\step[1]\id\\
\hh\step[-.5]\id\hstep\hld\step[.5]\id\hstep\sw1\step[-.5]\d\step\d\hstep\id\step\id\\
\hh\step[-.5]\id\step[.5]\id\hstep\id\hstep\hlu\hstep\cd\hstep\hld\hstep\id\hstep\cd\\
\hh\step[-.5]\id\hstep\id\hstep\x\step[.5]\id\step\d\step[-.5]\dd\hstep\id\hstep\d\step[-.5]\dd\hstep\hstep\id\\
\hh\step[-.5]\id\step[.5]\id\hstep\id\hstep\hld\hstep\cu\hstep\hlu\hstep\id\hstep\cu\\
\hh\step[-.5]\id\hstep\hlu\step[.5]\id\hstep\nw1\step[-.5]\dd\step\dd\hstep\id\step\id\\
\hh\step[-.5]\cu\step[.5]\cu\hstep\cu\ne1\step[1]\id\\
\hh\hstep\step[-.5]\id\hstep\hstep\hstep\id\step\ne1\step[-.5]\d\step[2]\id\\
\hh\hstep\step[-.5]\id\hstep\step\cu\step\id\step[2]\id\\
\object{\sstyle B}\step[2]\object{\sstyle H}\step[1.5]\object{\sstyle B}\step[2]\object{\sstyle H}\\
\end{tangle}
\hstretch 100 \vstretch 100 \enspace=\enspace
\begin{tangle}
\step[.75]\object{\sstyle B}\step[2.5]\object{\sstyle H}\step[2.5]
\object{\sstyle B}\step[2.5]\object{\sstyle H}\\
\hh{\step[-.5]\hstr{200}\cd}\step{\hstr{150}\cd}\step{\hstr{150}\cd}\step {\hstr{150}\cd}\\
\hh\step[-.5]\id\step{\hstr{200}\hld\hstep\id\hstep}\hstep\id\step\id\hstep{\hstr{200}\hld\hstep\id\hstep}\hstep\id\\
\hh\step[-.5]\id\step\id\step\se1\step[-.5]\dd\step\hstep\id\step\id\hstep\id\step\x\step\hstep\id\\
\hh\step[-.5]\id\hstep\cd\hstep\cd\d\hstep\hstep\x\hstep\cu\step\id\step\hstep\id\\
\hh\step[-.5]\id\hstep\id\step\d\step[-.5]\dd\step\id\hstep\x\step\x\step\hstep\id\hstep\step\id\\
\hh\step[-.5]\id\hstep\cu\hstep\cu\dd\hstep\hstep\x\hstep\cd\step\id\step\hstep\id\\
\hh\step[-.5]\id\hstep\hstep\id\step\ne1\step[-.5]\d\step\hstep\id\step
        \id\hstep\id\step\x\step\hstep\id\\
\hh\step[-.5]\id\hstep\hstep{\hstr{200}\hlu\hstep\id\hstep}\hstep\id\step
        \id\hstep{\hstr{200}\hlu\hstep\id\hstep}\hstep\id\\
\hh{\step[-.5]\hstr{200}\cu}\step{\hstr{150}\cu}\step{\hstr{150}\cu}\step
        {\hstr{150}\cu}\\
\step[.75]\object{\sstyle B}\step[2.5]\object{\sstyle H}\step[2.5]
\object{\sstyle B}\step[2.5]\object{\sstyle H}\\
\end{tangle}
\enspace=\enspace
\begin{tangle}
\step[.75]\object{\sstyle B}\step[2.5]\object{\sstyle H}\step[2.5]
\object{\sstyle B}\step[2.5]\object{\sstyle H}\\
\hh{\hstr{150}\cd}\step{\hstr{150}\cd}\step{\hstr{150}\cd}\step {\hstr{150}\cd}\\
\hh\id\hstep{\hstr{200}\hld\hstep\id\hstep}\hstep\id\step\id\hstep{\hstr{200}\hld\hstep\id\hstep}\hstep\id\\
\hh\id\hstep\id\step\x\step\hstep\id\step\id\hstep\id\step\x\step\hstep\id\\
\hh\id\hstep\cu\step\d\hstep\hstep\x\hstep\cu\step\id\step\hstep\id\\
\hh\id\step\id\Step\x\step\x\step\hstep\id\hstep\step\id\\
\hh\id\hstep\cd\step\dd\hstep\hstep\x\hstep\cd\step\id\step\hstep\id\\
\hh\id\hstep\id\step\x\step\hstep\id\step
        \id\hstep\id\step\x\step\hstep\id\\
\hh\id\hstep{\hstr{200}\hlu\hstep\id\hstep}\hstep\id\step
        \id\hstep{\hstr{200}\hlu\hstep\id\hstep}\hstep\id\\
\hh{\hstr{150}\cu}\step{\hstr{150}\cu}\step{\hstr{150}\cu}\step
        {\hstr{150}\cu}\\
\step[.75]\object{\sstyle B}\step[2.5]\object{\sstyle H}\step[2.5]
\object{\sstyle B}\step[2.5]\object{\sstyle H}
\end{tangle}
\]
\caption{Proof of bialgebra structure of $B\star H$}
\label{fig:biproduct}
\end{figure}
%**********************************************************************
\hstretch 50 \vstretch 50
\[
\begin{tangle}
\step[3.5]\object{\sstyle B}\Step\step[1.5]\object{\sstyle H}\\
\step[2]{\hstr{75}\cd}\step\cd\\
\step[2]\id\step\Ld\step\id\Step\id\\
\Ld\id\step\id\Step\hx\Step\id\\
\id\step[2]\id\step\cu\step\id\Step\id\\
\id\Step\x\Step\id\Step\id\\
\cu\Step\id\Step\id\Step\id\\
\obox 2{\sstyle {S_H}}\step\obox 2{\sstyle {S_B}}\step\id\Step\id\\
\cd\Step\id\Step\id\Step\id\\
\id\Step\x\Step\id\Step\id\\
\id\step[2]\id\step\cd\step\id\Step\id\\
\Lu\id\step\id\Step\hx\Step\id\\
\step[2]\id\step\Lu\step\id\Step\id\\
\step[2]{\hstr{75}\cu}\step\cu\\
\step[3.5]\object{\sstyle B}\step[1.5]\Step\object{\sstyle H}\\
\end{tangle}
\step=\step
\begin{tangle}
\step[3]\object{\sstyle B}\Step\Step\object{\sstyle H}\\
\step\Ld\Step\step\cd\\
\step\id\step\cd\Step\id\Step\id\\
\step\id\step\d\step\x\Step\id\\
\step\id\Step\hx\Step\id\Step\id\\
\step\cu\step\id\Step\id\Step\id\\
\step\obox 2{\sstyle {S_H}}\obox 2{\sstyle {S_B}}\step\id\Step\id\\
\step\cd\step\id\Step\id\Step\id\\
\cd\step\hx\Step\id\Step\id\\
\id\Step\hx\step\d\step\id\Step\id\\
\Lu\step\id\Step\hx\Step\id\\
\step[2]\id\step\Lu\step\cu\\
\step[2]{\hstr{75}\cu}\step[2]\id\\
\step[3.3]\object{\sstyle B}\step[1.5]\Step\object{\sstyle H}\\
\end{tangle}
\step=\step
\begin{tangle}
\step[2]\object{\sstyle B}\step\Step\object{\sstyle H}\\
\Ld\Step\cd\\
\id\step\cd\step\id\Step\id\\
\id\step\d\step\hx\Step\id\\
\id\Step\hx\step\id\Step\id\\
\cu\step\id\step\id\Step\id\\
\obox 2{\sstyle {S_H}}\obox 2{\sstyle {S_B}}\id\Step\id\\
\cd\step\id\step\id\Step\id\\
\id\Step\hx\step\id\Step\id\\
\id\step\dd\step\hx\Step\id\\
\id\step\cu\step\cu\\
\Lu\Step\step\id\\
\step[2]\object{\sstyle B}\step\Step\object{\sstyle H}\\
\end{tangle}
\step=\step
\begin{tangle}
\step[2]\object{\sstyle B}\step\Step\object{\sstyle H}\\
\Ld\Step\cd\\
\id\step[2]\QQ \ve\step\dd\step[2]\id\\
\id\Step\dd\step\Step\id\\
\cu\step\step\Step\id\\
\obox 2{\sstyle {S_H}}\step[2]\Step\id\\
\cd\step\step\Step\id\\
\id\Step\d\step\Step\id\\
\id\step[2]\Q \eta\step\d\step[2]\id\\
\Lu\Step\cu\\
\step[2]\object{\sstyle B}\step\Step\object{\sstyle H}\\
\end{tangle}
\step=\step
\begin{tangle}
\object{\sstyle B}\step[2]\Step\object{\sstyle H}\\
\QQ \ve\step[3]\cd\\
\Q \eta\Step\dd\step\step\id\\
\cu\step\Step\id\\
\obox 2{\sstyle {S_H}}\step\Step\id\\
\cd\step\Step\id\\
\QQ \ve\Step\d\step\step\id\\
\Q \eta\step[3]\cu\\
\object{\sstyle B}\step[2]\Step\object{\sstyle H}\\
\end{tangle}
\step=\step
\begin{tangle}
\object{\sstyle B}\step[3]\object{\sstyle H}\\
\id\step[2]\cd\\
\QQ \ve \Step\id \Step\id\\
\step\obox 2{\sstyle {S_H}}\step\id\\
\Q \eta \Step\id \Step\id\\
\id\step[2]\cu\\
\object{\sstyle B}\step[3]\object{\sstyle H}\\
\end{tangle}
\step=\step
\begin{tangle}
\object{\sstyle B}\Step\object{\sstyle H}\\
\id\Step\id\\
\QQ \ve \Step\QQ \ve \\
\Q \eta\Step\Q \eta\\
\id\Step\id\\
\object{\sstyle B}\Step\object{\sstyle H}\\
\end{tangle}
\]
%%%%%%%%%%%%%%%%%%%%%%%%%%%%%%%%%%%%%%%%%%%%%%%%%%%%%%%%%%%%%%%%%%%%%%%%%%%
\begin{figure}[!ht]
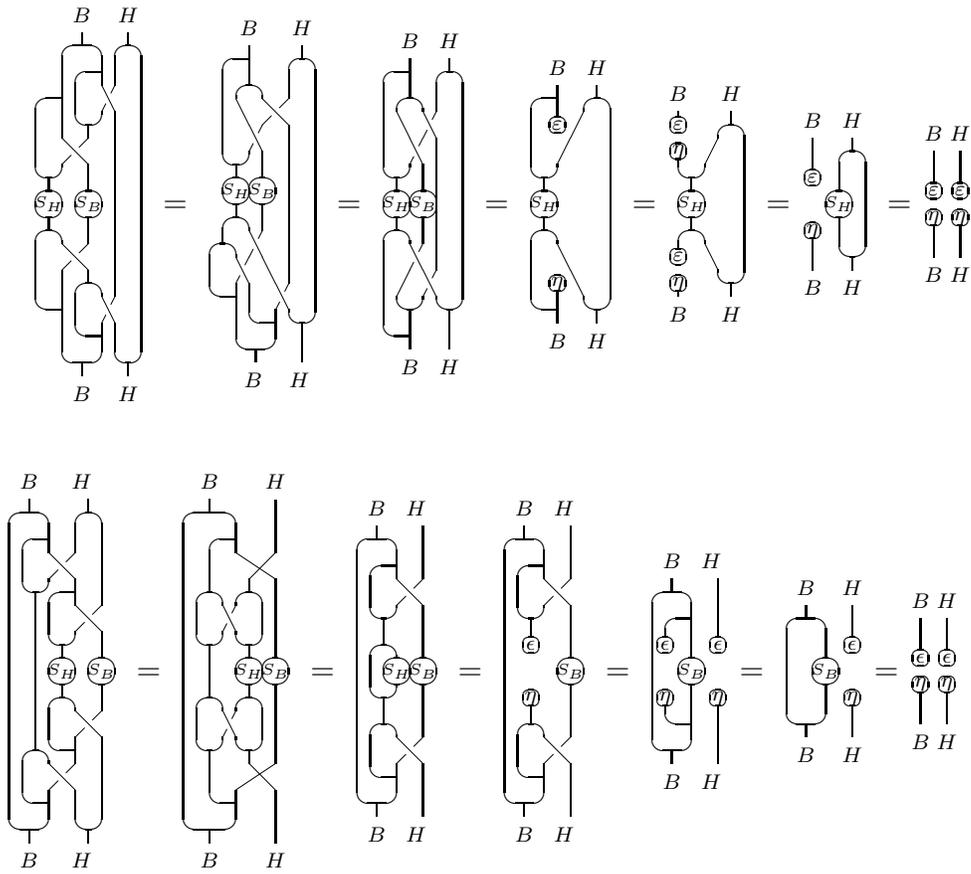

\[
\begin{tangle}
\step[1.5]\object{\sstyle B}\Step\step[2]\object{\sstyle H}\\
{\hstr{75}\cd}\Step\cd\\
\id\step\Ld\Step\id\Step\id\\
\id\step\id\Step\x\Step\id\\
\id\step\cu\Ld\Step\id\\
\id\Step\id\step\id\Step\x\\
\id\Step\id\step\cu\Step\id\\
\id\Step\id\step\obox 2{\sstyle {S_H}}\step\obox 2{\sstyle {S_B}}\\
\id\Step\id\step\cd\Step\id\\
\id\Step\id\step\id\Step\x\\
\id\step\cd\Lu\Step\id\\
\id\step\id\Step\x\Step\id\\
\id\step\Lu\Step\id\Step\id\\
{\hstr{75}\cu}\Step\cu\\
\step[1.5]\object{\sstyle B}\Step\step[2]\object{\sstyle H}\\
\end{tangle}
\step=\step
\begin{tangle}
\step[2]\object{\sstyle B}\Step\step[3]\object{\sstyle H}\\
\Cd\Step\step\id\\
\id\step[2]\Ld\step[3]\id\\
\id\Step\id\step\Step\step\ne2\step[-3]\nw3\\
\id\step\cd\step\cd\step\id\\
\id\step\id\Step\hx\Step\id\step\id\\
\id\step\cu\step\cu\step\id\\
\id\Step\id\Step\obox 2{\sstyle {S_H}}\obox 2{\sstyle {S_B}}\\
\id\step\cd\step\cd\step\id\\
\id\step\id\Step\hx\Step\id\step\id\\
\id\step\cu\step\cu\step\id\\
\id\Step\id\step\Step\step\se2\step[-3]\sw3\\
\id\step[2]\Lu\step[3]\id\\
\Cu\Step\step\id\\
\step[2]\object{\sstyle B}\Step\step[3]\object{\sstyle H}\\
\end{tangle}
\step=\step
\begin{tangle}
\step[1.5]\object{\sstyle B}\step\Step\object{\sstyle H}\\
{\hstr{75}\cd}\Step\id\\
\id\step\Ld\Step\id\\
\id\step\id\Step\x\\
\id\step\cu\Step\id\\
\id\step\cd\Step\id\\
\id\step\id\step\obox 2{\sstyle {S_H}}\obox 2{\sstyle {S_B}}\\
\id\step\cu\Step\id\\
\id\step\cd\Step\id\\
\id\step\id\Step\x\\
\id\step\Lu\Step\id\\
{\hstr{75}\cu}\Step\id\\
\step[1.5]\object{\sstyle B}\step\Step\object{\sstyle H}\\
\end{tangle}
\step=\step
\begin{tangle}
\step[1.5]\object{\sstyle B}\step\Step\object{\sstyle H}\\
{\hstr{75}\cd}\Step\id\\
\id\step\Ld\Step\id\\
\id\step\id\Step\x\\
\id\step\cu\Step\id\\
\id\step\step\QQ \epsilon\step\Step\id\\
\id\step\step[3]\obox 2{\sstyle {S_B}}\\
\id\step\step\Q \eta \step\Step\id\\
\id\step\cd\Step\id\\
\id\step\id\Step\x\\
\id\step\Lu\Step\id\\
{\hstr{75}\cu}\Step\id\\
\step[1.5]\object{\sstyle B}\step\Step\object{\sstyle H}\\
\end{tangle}
\step=\step
\begin{tangle}
\step[1.5]\object{\sstyle B}\step\Step\object{\sstyle H}\\
{\hstr{75}\cd}\Step\id\\
\id\step\Ld\Step\id\\
\id\step\QQ \epsilon\step[2]\id\step[2]\QQ \epsilon\\
\id\step[2]\obox 2{\sstyle {S_B}}\\
\id\step\Q \eta \step[2]\id\step[2]\Q \eta\\
\id\step\Lu\Step\id\\
{\hstr{75}\cu}\Step\id\\
\step[1.5]\object{\sstyle B}\step\Step\object{\sstyle H}\\
\end{tangle}
\step=\step
\begin{tangle}
\step[1.5]\object{\sstyle B}\step[1.5]\Step\object{\sstyle H}\\
{\hstr{75}\cd}\Step\id\\
\id\step\step[2]\id\Step\QQ \epsilon\\
\id\step[2]\obox 2{\sstyle {S_B}}\\
\id\step\step[2]\id\Step\Q \eta\\
{\hstr{75}\cu}\Step\id\\
\step[1.5]\object{\sstyle B}\step[1.5]\Step\object{\sstyle H}\\
\end{tangle}
\step=\step
\begin{tangle}
\object{\sstyle B}\Step\object{\sstyle H}\\
\id\Step\id\\
\QQ \epsilon \Step\QQ \epsilon \\
\Q \eta\Step\Q \eta\\
\id\Step\id\\
\object{\sstyle B}\Step\object{\sstyle H}\\
\end{tangle}\ \ \ .
\]
\caption{Proof of antipode of $B\star H$} \label{fig:antipode}
\end{figure}

\vskip0.1cm

\vskip0.1cm

College of Mathematics, Henan Normal University, Xinxiang 453007, China

E-mail address: zhangtao@htu.cn


\providecommand{\bysame}{\leavevmode\hbox
to3em{\hrulefill}\thinspace}



\newpage

\begin{thebibliography}{Ra85}

\bibitem{AS02} N. Andruskiewitsch and H.-J. Schneider, Lifting of Quantum Linear Spaces and
Pointed Hopf Algebras of order $p^3$, J. Algebra 209 (1998), 658--691.


\bibitem{BD99} Y.~Bespalov, B.~Drabant. Cross product bialgebras - Part I, J. Algebra, 219(1999), 466--505.


\bibitem{Ra85}
D.~E.~ Radford, The Structure of Hopf Algebras with a Projection,
J.~Algebra 92(1985), 322--374.


\bibitem{Ma94} S. Majid,
Cross Products by Braided Groups and Bosonization, J. Algebra 163
(1994), 165--190.

\bibitem{Mo93}S. Montgomery, Hopf Algebras and Their Actions on Rings, CBMS 82,
American Math Society, Providence, RI, 1993.


\bibitem{Ye90} D. N. Yetter, Quantum groups and representations of monoidal
categorics, Math. Proc. Cambridge Phil. Soc., 108 (1990), 261-290.


\bibitem{ZC01} S.C. Zhang, H.X. Chen, The double bicrossproducts in
braided tensor categories, Comm. in Algebra 29(2001), 31--66.

%\bibitem{Zh99} Zhang S C. Braided Hopf Algebras. Changsha: Hunan Normal University Press, 1999 arXiv:math.RA/0511251.

\end{thebibliography}
\end{document}